\documentclass[final,3p,11pt,nopreprintline]{elsarticle}

\usepackage{amsmath,amssymb,amsthm,mathtools}
\usepackage{lmodern}
\usepackage{microtype}
\usepackage[hidelinks]{hyperref}

\journal{Topology and its Applications}
\hypersetup{
  pdftitle={Four problems about the irreducible determined T0-topologies},
  pdfauthor={Chongyun Zhao and Xiaoquan Xu}
}

\newtheorem{theorem}{Theorem}[section]
\newtheorem{proposition}[theorem]{Proposition}
\newtheorem{lemma}[theorem]{Lemma}
\newtheorem{corollary}[theorem]{Corollary}
\theoremstyle{definition}
\newtheorem{definition}[theorem]{Definition}
\newtheorem{example}[theorem]{Example}
\theoremstyle{remark}
\newtheorem{remark}[theorem]{Remark}

\DeclareMathOperator{\Irr}{Irr}
\DeclareMathOperator{\cl}{cl}
\DeclareMathOperator{\rd}{rd}
\DeclareMathOperator{\Rd}{Rd}
\DeclareMathOperator{\Hd}{Hd}
\DeclareMathOperator{\Fin}{Fin}
\newcommand{\RDset}{\mathsf{RD}}
\newcommand{\WDset}{\mathsf{WD}}
\newcommand{\Pow}{\mathcal P}
\newcommand{\Izero}{\mathcal I_0}
\newcommand{\Nat}{\mathbb N}
\newcommand{\up}{\mathord{\uparrow}}
\newcommand{\dn}{\mathord{\downarrow}}

\begin{document}
\begin{frontmatter}

\title{Four problems about the irreducible determined $T_0$-topologies}

\author[scu]{Chongyun Zhao}
\author[nit]{Xiaoquan Xu\corref{cor}}
\ead{xiqxu2002@163.com}
\cortext[cor]{Corresponding author.}
\address[scu]{School of Mathematics, Sichuan University, Chengdu 610065, China}
\address[nit]{School of Computer Information Engineering, Nanchang Institute of Technology, Nanchang 330044, China}

\begin{abstract}
We answer four questions on irreducible determined $T_0$-topologies and a related question on monotone determined topologies posed by Xu, Yang and Luo (Topology Appl.\ 382 (2026), 109762).
For every $T_0$ topology $\tau$, the family $\tau_r$ used to define its irreducible-derived topology is already a topology and preserves all nonempty irreducible subsets and their closures.
Consequently, the least irreducible determined refinement is obtained in one step:
$\tau_r=\rd(\tau)=\Rd(\tau)$.
For an arbitrary poset $P$, we describe the weakest order-compatible irreducible determined topology $\rho(P)$ using Dedekind--MacNeille cuts and a finite-cover condition involving principal lower sets.
The Scott space of the binary-prefix domain provides a common counterexample to saturated heredity for irreducible determined spaces and monotone determined spaces: its maximal-point subspace is Cantor space.
We also describe $\rho(P)$-closure by iterating an irreducible-cut operator, characterize Scott continuity of this operator, and give an example requiring two applications of the cut operator.
A finite-decomposition condition extends the topology and closure-preservation results to general subset systems.
For the Rudin and well-filtered determined systems, we exhibit a common $T_1$ example whose first and second derived topologies are the co-countable and discrete topologies, respectively.
Thus the first derived topology need not be determined by the same subset system.
\end{abstract}


\end{frontmatter}

\section{Introduction}\label{sec:introduction}

The Scott topology and the specialization order connect partially ordered sets with $T_0$ spaces and provide a setting for extending domain-theoretic ideas beyond ordered structures; see \cite{CLD}.
One such extension uses \emph{monotone determined spaces}, or MD-spaces, whose open sets are determined by the closures of directed subsets.
This notion originates in the work of Ern\'e \cite{Erne2009}; its derived topologies were studied by Luo and Xu \cite{LuoXu2017}.
Yu and Kou independently introduced \emph{directed spaces} through convergence with respect to the specialization order \cite{YuKou2015}.
Directed spaces and MD-spaces coincide, and include the Scott spaces of all posets.

Every nonempty directed subset of a $T_0$ space is irreducible, so it is natural to replace directed subsets by irreducible subsets in topology-generating constructions.
Zhao and Ho introduced the $\mathrm{SI}$-topology and related it to $k$-bounded sobriety \cite{ZhaoHo}.
Cuts provide a way to formulate analogous constructions without requiring the existence of suprema, as in Ern\'e's work on Scott topology for posets \cite{Erne1981} and Luo and Xu's study of $\mathrm{SI}_2$-continuous spaces \cite{LuoXu2019}.
The $\mathrm{SI}$ and $\mathrm{SI}_2$ constructions produce coarser topologies.
The irreducible determined construction considered in this paper instead refines the given topology.

Xu, Yang and Luo \cite{XuYangLuo} introduced \emph{irreducible determined spaces}, or RD-spaces, and the irreducible-derived topology $\rd(\tau)$ of a $T_0$ topology $\tau$.
They proved that every MD-space is an RD-space and that there exists a least RD topology $\Rd(\tau)$ containing $\tau$.
They also obtained a transfinite description of $\Rd(\tau)$ and established that RD-spaces are closed under open subspaces, closed subspaces and retracts.
Their work leaves open several questions concerning irreducible subsets, the refinement operation, order-compatible topologies and saturated subspaces.

This paper answers the four RD questions in \cite[Questions~6.1--6.4]{XuYangLuo}, together with the related MD question in \cite[Question~3.14]{XuYangLuo}.
The contributions are as follows.
\begin{enumerate}
\item We prove that the family $\tau_r$ defining $\rd(\tau)$ is already a topology.
It has precisely the same nonempty irreducible subsets as $\tau$, and preserves the closure of each such subset.
It follows that
\[
 \rd(\rd(\tau))=\rd(\tau)=\Rd(\tau).
\]
Thus the least RD refinement is obtained in one step, with no additional assumptions on the $T_0$ space.
These results answer Questions~6.1 and~6.2; see Theorems~\ref{thm:preservation} and~\ref{thm:one-step}.
\item For a poset $P$, we describe the weakest RD topology having the prescribed specialization order.
The description is entirely order-theoretic: it uses cuts of nonempty subsets satisfying a finite-cover condition by principal lower sets.
The resulting topology $\rho(P)$ satisfies
\[
 \rho(P)=\rd(\upsilon(P))=\Rd(\upsilon(P)),
\]
where $\upsilon(P)$ is the upper topology.
This answers Question~6.3; see Theorem~\ref{thm:rho}.
\item The maximal elements of the binary-prefix domain form a saturated subspace of its Scott space that is homeomorphic to Cantor space.
The ambient space is both MD and RD, whereas the subspace has neither property.
This gives negative answers to Question~6.4 and the additional MD question, Question~3.14; see Theorem~\ref{thm:saturated}.
\end{enumerate}

We then study the irreducible-cut operator associated with $\rho(P)$.
Its fixed points are exactly the $\rho(P)$-closed sets, and its transfinite iteration computes closure.
We characterize its Scott continuity and show, by an example on an algebraic complete lattice, that the operator itself need not be idempotent.
This distinguishes the one-step construction of an RD topology from the computation of closures of arbitrary subsets.

We also extend the finite-partition argument to the subset-system derived topologies of Xu, Huang and Yang \cite{XuHuangYang2026}.
A finite-decomposition condition ensures that the defining family $\tau_H$ is a topology and preserves the closures of the original $H$-sets.
An additional preservation condition implies that the first derived topology is $H$-determined.
The decomposition condition holds for the singleton, chain, directed, Rudin, well-filtered determined and irreducible subset systems.
For the Rudin and well-filtered determined systems, we construct a common $T_1$ example whose derived topologies stabilize after exactly two steps.
These results give a sufficient condition and two non-idempotent standard cases for the characterization problem in \cite[Question~7.1]{XuHuangYang2026}; see Section~\ref{sec:subset-systems}.

Section~\ref{sec:preliminaries} records the notation and a finite-partition lemma.
Section~\ref{sec:refinement} proves the one-step refinement theorem.
Section~\ref{sec:rho} gives the order-theoretic description of $\rho(P)$.
Section~\ref{sec:saturated} treats the binary-prefix domain and saturated subspaces.
Section~\ref{sec:operator} studies the cut operator.
Section~\ref{sec:subset-systems} treats general subset systems and the Rudin and well-filtered determined examples, and Section~\ref{sec:conclusion} concludes the paper.

\section{Preliminaries}\label{sec:preliminaries}

All directed and irreducible subsets in this paper are assumed to be nonempty.
We use $\Nat=\{0,1,2,\ldots\}$.
For a poset $P$ and $A\subseteq P$, write
\[
\begin{aligned}
 \dn A&=\{x\in P:x\leq a\text{ for some }a\in A\},
 &\up A&=\{x\in P:a\leq x\text{ for some }a\in A\},\\
 A^u&=\{x\in P:a\leq x\text{ for all }a\in A\},
 &A^l&=\{x\in P:x\leq a\text{ for all }a\in A\}.
\end{aligned}
\]
We abbreviate $\dn\{x\}$ and $\up\{x\}$ to $\dn x$ and $\up x$.
The \emph{Dedekind--MacNeille cut} of $A$ is
\[
 A^\delta=(A^u)^l=\bigcap_{p\in A^u}\dn p.
\]
An empty intersection is understood to be $P$; in particular, $A^\delta=P$ if $A$ has no upper bound.
We always have $A\subseteq A^\delta$, and $A^\delta=\dn\sup A$ whenever $\sup A$ exists.
The set of all subsets of $P$ is denoted by $\Pow(P)$, and $\Fin(P)$ denotes its finite subsets, including $\varnothing$.

A subset $D\subseteq P$ is \emph{directed} if it is nonempty and every two of its elements have an upper bound in $D$.
A \emph{dcpo} is a poset in which every directed subset has a supremum.
An element $k$ of a dcpo is \emph{compact} if $k\leq\sup D$ for a directed set $D$ implies $k\leq d$ for some $d\in D$.
A dcpo is \emph{algebraic} if, for every $x$, the compact elements below $x$ form a directed set with supremum $x$.

The \emph{upper topology} $\upsilon(P)$ has subbase
$\{P\setminus\dn p:p\in P\}$ and base
$\{P\setminus\dn F:F\in\Fin(P)\}$.
The \emph{Alexandrov topology} $\alpha(P)$ consists of all upper subsets of $P$.
The \emph{Scott topology} $\sigma(P)$ consists of the upper sets $U$ such that
\[
 \sup D\in U~\Longrightarrow~ D\cap U\neq\varnothing
\]
for every directed $D$ whose supremum exists.
We also use the cut version of the Scott topology, denoted by $\sigma_2(P)$:
\begin{equation}\label{eq:sigma2}
 \sigma_2(P)=\{U=\up U:
   U\cap D^\delta\neq\varnothing\Longrightarrow U\cap D\neq\varnothing
   \text{ for every directed }D\subseteq P\}.
\end{equation}
This is the Scott--Ern\'e, or weak Scott, topology; see \cite{Erne1981,LuoXu2019}.
For dcpos, $\sigma_2(P)=\sigma(P)$.

For a $T_0$ space $X=(X,\tau)$, we use the specialization order
\[
 x\leq_X y~\Longleftrightarrow~ x\in\cl_\tau\{y\}.
\]
A topology on a poset $P$ is \emph{order-compatible} if its specialization order is the given order on $P$; equivalently,
\begin{equation}\label{eq:compatible}
 \upsilon(P)\subseteq\tau\subseteq\alpha(P).
\end{equation}
Every open set is upper for the specialization order.
The saturated subsets, that is, intersections of open sets, are precisely the upper subsets.

A nonempty subset $A$ of a space $X$ is \emph{irreducible} if
$A\subseteq F_1\cup F_2$, with $F_1,F_2$ closed, implies $A\subseteq F_1$ or $A\subseteq F_2$.
Equivalently, whenever two open sets meet $A$, their intersection meets $A$.
We write $\Irr(X,\tau)$ for the family of all nonempty irreducible subsets, not only the closed ones.
A subset is irreducible if and only if its closure is irreducible.
Every directed subset of the specialization order is irreducible, since open sets are upper.

The following elementary observation is used repeatedly.

\begin{lemma}[Finite-partition lemma]\label{lem:partition}
Let $A$ be an irreducible subset of a topological space $X$.
If $A\subseteq C_1\cup\cdots\cup C_n$, where $n\geq1$ and the $C_i$ are arbitrary subsets of $X$, then for some $i$ the set $B=A\cap C_i$ is nonempty and irreducible, and $\cl B=\cl A$.
\end{lemma}

\begin{proof}
The equality
$\cl A=\bigcup_{i=1}^n\cl(A\cap C_i)$ and irreducibility imply
$\cl A=\cl(A\cap C_i)$ for some $i$.
The common closure is nonempty and irreducible, so $A\cap C_i$ has the required properties.
\end{proof}

For completeness, put
\[
 \tau_d=\{U\subseteq X:
  U\cap\cl_\tau D\neq\varnothing\Longrightarrow U\cap D\neq\varnothing
  \text{ for every directed }D\subseteq X\}.
\]
The space $(X,\tau)$ is \emph{monotone determined} if $\tau=\tau_d$.
Singleton tests imply that every member of $\tau_d$ is upper.
Every Scott space $(P,\sigma(P))$ is MD: if $U\in\sigma(P)_d$ and a directed set $D$ has supremum in $U$, then
$\sup D\in\cl_{\sigma(P)}D$, so $D$ meets $U$ and $U$ is Scott open.
The reverse inclusion $\sigma(P)\subseteq\sigma(P)_d$ follows from the definition of closure.

\section{One-step irreducible determination}\label{sec:refinement}

We recall the construction of Xu, Yang and Luo \cite{XuYangLuo}.
For a $T_0$ topology $\tau$ on $X$, let
\begin{equation}\label{eq:tau-r}
 \tau_r=\{U\subseteq X:
  U\cap\cl_\tau A\neq\varnothing\Longrightarrow U\cap A\neq\varnothing
  \text{ for every }A\in\Irr(X,\tau)\}.
\end{equation}
The topology generated by $\tau_r$ is denoted by $\rd(\tau)$.
The space $(X,\tau)$, or the topology $\tau$, is \emph{irreducible determined} if $\rd(\tau)=\tau$; we abbreviate this to RD.
Since $\tau\subseteq\tau_r$, the RD condition is equivalent to $\tau_r=\tau$.
Moreover, $\tau\subseteq\tau_r\subseteq\tau_d$, so every MD-space is RD.

The existence of a least RD topology containing $\tau$, denoted by $\Rd(\tau)$, was established in \cite{XuYangLuo} by a transfinite construction.
We show that it is already $\tau_r$.

\begin{proposition}\label{prop:tau-r-topology}
For every $T_0$ topology $\tau$ on $X$, the family $\tau_r$ is a topology.
It contains $\tau$ and has the same specialization order as $\tau$.
Its closed sets are exactly the members of
\begin{equation}\label{eq:closed-family}
 \mathcal F_\tau=\{C\subseteq X:
  A\in\Irr(X,\tau),\ A\subseteq C
  \Longrightarrow\cl_\tau A\subseteq C\}.
\end{equation}
In particular, $\rd(\tau)=\tau_r$.
\end{proposition}

\begin{proof}
Taking complements in \eqref{eq:tau-r} gives
$U\in\tau_r$ if and only if $X\setminus U\in\mathcal F_\tau$.
Both $\varnothing$ and $X$ belong to $\mathcal F_\tau$, and this family is closed under arbitrary intersections.
To verify closure under binary unions, let $C_1,C_2\in\mathcal F_\tau$ and let $A\in\Irr(X,\tau)$ be contained in $C_1\cup C_2$.
By Lemma~\ref{lem:partition}, for some $i\in\{1,2\}$ the set $B=A\cap C_i$ is irreducible and $\cl_\tau B=\cl_\tau A$.
Hence
\[
 \cl_\tau A=\cl_\tau B\subseteq C_i\subseteq C_1\cup C_2.
\]
Thus $\mathcal F_\tau$ is the closed-set family of a topology, namely $\tau_r$.

Every $\tau$-open set satisfies \eqref{eq:tau-r}, so $\tau\subseteq\tau_r$.
For $U\in\tau_r$, applying \eqref{eq:tau-r} to singletons shows that $U$ is upper for the specialization order of $\tau$.
Consequently, $\tau\subseteq\tau_r\subseteq\alpha(X,\leq_\tau)$, and the specialization orders agree.
\end{proof}

\begin{theorem}[Preservation of irreducibles and their closures]\label{thm:preservation}
Let $\tau$ be a $T_0$ topology on $X$.
Then
\begin{equation}\label{eq:preserve-irr}
 \Irr(X,\rd(\tau))=\Irr(X,\tau).
\end{equation}
For every $A\in\Irr(X,\tau)$,
\begin{equation}\label{eq:preserve-closure}
 \cl_{\rd(\tau)}A=\cl_\tau A.
\end{equation}
\end{theorem}

\begin{proof}
Since $\tau\subseteq\rd(\tau)$, every $\rd(\tau)$-irreducible subset is $\tau$-irreducible.
Conversely, let $A\in\Irr(X,\tau)$, and suppose that $A\subseteq C_1\cup C_2$ with $C_1,C_2$ closed in $\rd(\tau)$.
Lemma~\ref{lem:partition} gives an irreducible set $B=A\cap C_i$ such that $\cl_\tau B=\cl_\tau A$ for some $i$.
By \eqref{eq:closed-family},
\[
 A\subseteq\cl_\tau A=\cl_\tau B\subseteq C_i.
\]
Thus $A$ is $\rd(\tau)$-irreducible, proving \eqref{eq:preserve-irr}.

For the closure identity, every $\rd(\tau)$-closed set containing $A$ contains $\cl_\tau A$, by \eqref{eq:closed-family}.
Hence $\cl_\tau A\subseteq\cl_{\rd(\tau)}A$.
The reverse inclusion follows from $\tau\subseteq\rd(\tau)$.
\end{proof}

\begin{theorem}[One-step refinement]\label{thm:one-step}
For every $T_0$ topology $\tau$,
\begin{equation}\label{eq:one-step}
 \rd(\rd(\tau))=\rd(\tau)=\Rd(\tau)=\tau_r.
\end{equation}
Moreover, if $\tau\subseteq\theta$ are $T_0$ topologies on the same set, then
$\rd(\tau)\subseteq\rd(\theta)$.
\end{theorem}

\begin{proof}
The defining test in \eqref{eq:tau-r} involves only the irreducible subsets and their closures.
These are unchanged by Theorem~\ref{thm:preservation}, so
$\rd(\rd(\tau))=\rd(\tau)$.

If $\tau\subseteq\theta$, then
\[
 \Irr(X,\theta)\subseteq\Irr(X,\tau),
 ~ \cl_\theta A\subseteq\cl_\tau A.
\]
It follows directly from \eqref{eq:tau-r} that $\tau_r\subseteq\theta_r$.
Thus $\rd$ is monotone with respect to inclusion of topologies.
If $\theta$ is RD and contains $\tau$, then
$\rd(\tau)\subseteq\rd(\theta)=\theta$.
Since $\rd(\tau)$ is itself RD and contains $\tau$, it is the least such topology, namely $\Rd(\tau)$.
\end{proof}

Theorems~\ref{thm:preservation} and~\ref{thm:one-step} answer \cite[Questions~6.1 and~6.2]{XuYangLuo} without additional hypotheses.
The operation $\rd$ is therefore an extensive, monotone and idempotent operation on the $T_0$ topologies on a fixed set.

\section{The weakest order-compatible RD topology}\label{sec:rho}

We now express $\rd(\upsilon(P))$ using only the order on $P$.

\begin{definition}\label{def:i0}
Let $\Izero(P)$ be the family of all nonempty $A\subseteq P$ satisfying
\begin{equation}\label{eq:i0}
 A\subseteq\dn F,~ F\in\Fin(P)
 ~\Longrightarrow~
 A\subseteq\dn f\text{ for some }f\in F.
\end{equation}
\end{definition}

\begin{proposition}\label{prop:i0}
For every poset $P$,
\[
 \Izero(P)=\Irr(P,\upsilon(P)).
\]
If $A\in\Izero(P)$, then $\cl_{\upsilon(P)}A=A^\delta$.
In particular, every directed subset of $P$ belongs to $\Izero(P)$.
\end{proposition}

\begin{proof}
If $A$ is irreducible in $\upsilon(P)$, then every finite cover of $A$ by the closed sets $\dn f$ has a member containing $A$; hence \eqref{eq:i0} holds.
Conversely, assume \eqref{eq:i0}, and let $U,V$ be $\upsilon(P)$-open sets meeting $A$.
There are finite sets $F,G\subseteq P$ such that
$P\setminus\dn F\subseteq U$ and $P\setminus\dn G\subseteq V$, and both basic open sets meet $A$.
If their intersection did not meet $A$, we would have $A\subseteq\dn(F\cup G)$.
Condition~\eqref{eq:i0} would then imply $A\subseteq\dn h$ for some $h\in F\cup G$, contradicting one of the two nonempty intersections.
Thus $A$ is irreducible.

Every $\dn p$ with $p\in A^u$ is closed and contains $A$, so
$\cl_{\upsilon(P)}A\subseteq A^\delta$.
If $x\notin\cl_{\upsilon(P)}A$, choose a basic open neighbourhood
$P\setminus\dn F$ of $x$ disjoint from $A$.
Then $A\subseteq\dn F$, so $A\subseteq\dn f$ for some $f\in F$.
Since $x\notin\dn f$, we have $x\notin A^\delta$.
This proves the reverse inclusion.
Finally, directed subsets are irreducible in every order-compatible topology, in particular in $\upsilon(P)$.
\end{proof}

\begin{definition}\label{def:rho}
For a poset $P$, define
\begin{equation}\label{eq:rho}
 \rho(P)=\{U=\up U:
   U\cap A^\delta\neq\varnothing\Longrightarrow U\cap A\neq\varnothing
   \text{ for every }A\in\Izero(P)\}.
\end{equation}
We call $\rho(P)$ the \emph{upper irreducible cut topology} on $P$.
\end{definition}

\begin{samepage}
\begin{theorem}\label{thm:rho}
For every poset $P$, the following statements hold.
\begin{enumerate}[\rm (1)]
\item $\rho(P)=\rd(\upsilon(P))=\Rd(\upsilon(P))$.
It is the weakest order-compatible RD topology on $P$.
\item A subset $C\subseteq P$ is $\rho(P)$-closed if and only if
\begin{equation}\label{eq:rho-closed}
 A\in\Izero(P),\ A\subseteq C~\Longrightarrow~ A^\delta\subseteq C.
\end{equation}
\item $\Irr(P,\rho(P))=\Izero(P)$, and
$\cl_{\rho(P)}A=A^\delta$ for every $A\in\Izero(P)$.
\end{enumerate}
\end{theorem}
\end{samepage}

\begin{proof}
By Proposition~\ref{prop:i0}, the irreducible subsets and their closures in $\upsilon(P)$ are exactly the sets in $\Izero(P)$ and their cuts.
Thus the test in \eqref{eq:tau-r}, applied to $\upsilon(P)$, is the test in \eqref{eq:rho}.
The upper-set condition is automatic from singleton tests.
Proposition~\ref{prop:tau-r-topology} and Theorem~\ref{thm:one-step} give
$\rho(P)=\rd(\upsilon(P))=\Rd(\upsilon(P))$ and show that this topology is order-compatible.
Every order-compatible topology contains $\upsilon(P)$, so the minimality statement follows.

The closed-set description follows from \eqref{eq:closed-family} and Proposition~\ref{prop:i0}.
Using $A=\{c\}$ in \eqref{eq:rho-closed} gives $\dn c\subseteq C$ whenever $c\in C$.
The final assertion follows from Theorem~\ref{thm:preservation}.
\end{proof}

Theorem~\ref{thm:rho} answers \cite[Question~6.3]{XuYangLuo}.
It also shows that $\rho(P)$ can be defined without first specifying a topology on $P$.

\begin{proposition}\label{prop:chain}
For every poset $P$,
\begin{equation}\label{eq:chain}
 \upsilon(P)\subseteq\rho(P)\subseteq\sigma_2(P)
 \subseteq\sigma(P)\subseteq\alpha(P).
\end{equation}
If $P$ is a dcpo, then $\sigma_2(P)=\sigma(P)$.
\end{proposition}

\begin{proof}
The first inclusion follows from Theorem~\ref{thm:rho}.
Every directed set belongs to $\Izero(P)$, which gives $\rho(P)\subseteq\sigma_2(P)$.
If a directed set $D$ has a supremum, then $D^\delta=\dn\sup D$, proving $\sigma_2(P)\subseteq\sigma(P)$ and equality when $P$ is a dcpo.
The final inclusion is part of the definition of Scott openness.
\end{proof}

\begin{example}\label{ex:antichain}
Let $P$ be an infinite antichain.
Then $\upsilon(P)=\rho(P)$ is the cofinite topology, whereas
$\sigma_2(P)=\sigma(P)=\alpha(P)$ is discrete.
Indeed, $\Izero(P)$ consists of the singletons and the infinite subsets.
For an infinite $A$, we have $A^\delta=P$, while $\{p\}^\delta=\{p\}$.
Condition~\eqref{eq:rho} therefore says that every nonempty $\rho(P)$-open set meets every infinite subset of $P$, which is equivalent to having finite complement.
Every directed subset of an antichain is a singleton, so all the remaining topologies are discrete.
\end{example}

\begin{example}\label{ex:finite}
If $P$ is finite, all five topologies in \eqref{eq:chain} coincide.
To see that $\upsilon(P)=\alpha(P)$, let $U$ be an upper set.
Its complement is a finite lower set $F=\dn F$, so $U=P\setminus\dn F$ is upper-topology open.
The assertion follows from \eqref{eq:chain}.
\end{example}

\section{Saturated subspaces and the binary-prefix domain}\label{sec:saturated}

We use the same ordered set to illustrate the weakest RD topology and to answer the saturated-subspace questions.
The topology supplying the counterexample will be the Scott topology.

Let $2=\{0,1\}$ and
\[
 X=2^{<\omega}\cup2^\omega,
 ~ 2^{<\omega}=\bigcup_{n\in\Nat}2^n.
\]
Thus $X$ consists of all finite and infinite binary sequences, including the empty sequence.
Order $X$ by the prefix relation $s\sqsubseteq t$.
The finite prefix of length $n$ of an infinite sequence $\xi$ is denoted by $\xi|_n$.

\begin{proposition}\label{prop:binary-domain}
The prefix order makes $X$ an algebraic dcpo whose compact elements are exactly the finite sequences.
The sets $\up s$, $s\in2^{<\omega}$, form a base for $\sigma(X)$.
\end{proposition}

\begin{proof}
Two prefixes of a common sequence are comparable.
Consequently, every directed subset $D$ of $X$ is a chain, and its union is a finite or infinite sequence that is the supremum of $D$.
Hence $X$ is a dcpo.

If a finite sequence $s$ is a prefix of $\sup D$, then some member of $D$ has length at least the length of $s$ and extends $s$.
This proves that every finite sequence is compact; for the empty sequence, any member of the nonempty set $D$ suffices.
An infinite sequence $\xi$ is the supremum of the increasing chain
$\{\xi|_n:n\in\Nat\}$ and is not below any member of that chain, so it is not compact.
Every point is the supremum of its finite prefixes, proving algebraicity.

Each $\up s$ with $s$ finite is Scott open by compactness.
If $x$ belongs to a Scott-open set $U$, the directed set of its finite prefixes has supremum $x$, so one of these prefixes $s$ belongs to $U$.
Then $x\in\up s\subseteq U$, which proves the base assertion.
\end{proof}

\begin{proposition}\label{prop:binary-upper}
For the binary-prefix domain $X$,
\begin{equation}\label{eq:binary-upper}
 \rho(X)=\upsilon(X)
 =\{\varnothing\}\cup
   \{X\setminus\dn F:F\in\Fin(X)\}.
\end{equation}
Equivalently, the closed sets are $X$ and the finitely generated lower sets.
\end{proposition}

\begin{proof}
Let $C$ be a proper, nonempty $\rho(X)$-closed set.
If $C$ were not contained in any finite union of principal lower sets, it would satisfy \eqref{eq:i0} vacuously and would have no upper bound.
Thus $C\in\Izero(X)$ and $C^\delta=X$, contradicting \eqref{eq:rho-closed}.
Therefore $C\subseteq\dn F$ for some finite $F\subseteq X$.

For each $f\in F$ such that $C_f=C\cap\dn f$ is nonempty, the set $C_f$ is a chain with supremum $c_f=\bigcup C_f$.
Since directed sets belong to $\Izero(X)$, the closedness condition gives $c_f\in C$.
Also $c_f\sqsubseteq f$, and $C$ is lower, so $C_f=\dn c_f$.
Hence
\[
 C=\dn\{c_f:f\in F,\ C_f\neq\varnothing\}.
\]
The empty set is finitely generated as well.
Conversely, every finite union of principal lower sets is $\upsilon(X)$-closed and therefore $\rho(X)$-closed.
This identifies the closed sets and proves \eqref{eq:binary-upper}.
\end{proof}

\begin{theorem}\label{thm:saturated}
The properties of being an MD-space and being an RD-space are not hereditary with respect to saturated subspaces.
More precisely, $(X,\sigma(X))$ is both MD and RD, while its saturated subspace $M=2^\omega$ is homeomorphic to Cantor space and has neither property.
\end{theorem}

\begin{proof}
Every Scott space is MD, and every MD-space is RD.
Thus $(X,\sigma(X))$ has both properties.
The set $M$ consists of the maximal elements of $X$ and is upper, hence saturated.

For each finite binary sequence $s$, put
\[
 [s]=\{\xi\in2^\omega:s\sqsubseteq\xi\}
     =(\up s)\cap M.
\]
By Proposition~\ref{prop:binary-domain}, these cylinders form a base for the subspace topology on $M$.
They are precisely the standard basic open sets of the product space $2^\omega$ with $2$ discrete.
Thus $M$ is Cantor space; in particular, it is Hausdorff and not discrete.

Every nonempty irreducible subset of a Hausdorff space is a singleton, and its closure is itself.
Consequently, for a Hausdorff topology $\eta$, every subset satisfies the defining test for $\eta_r$, so $\eta_r$ is discrete.
A Hausdorff space is therefore RD if and only if it is discrete.
It follows that $M$ is not RD, and hence is not MD.
\end{proof}

This gives negative answers to \cite[Questions~3.14 and~6.4]{XuYangLuo} using one countably based algebraic dcpo.

\begin{remark}\label{rem:two-subspaces}
The choice of ambient topology in Theorem~\ref{thm:saturated} matters.
By Proposition~\ref{prop:binary-upper}, the topology induced by $\rho(X)=\upsilon(X)$ on $M$ is cofinite: a finitely generated lower set meets $M$ in a finite set, and every finite subset of $M$ arises in this way.
Since the order on $M$ is equality, Example~\ref{ex:antichain} shows that this cofinite topology is RD.
Thus the Scott subspace in Theorem~\ref{thm:saturated} and the upper-topology subspace on the same set are different.
\end{remark}

\section{The irreducible-cut operator and closure}\label{sec:operator}

The order-theoretic description of $\rho(P)$ gives an operator on $\Pow(P)$.
For $B\subseteq P$, define
\begin{equation}\label{eq:T}
 T_P(B)=\bigcup\{A^\delta:A\in\Izero(P),\ A\subseteq B\}.
\end{equation}
We call $T_P$ the \emph{irreducible-cut operator} and write $T$ when $P$ is clear.
Also write $c_P(B)=\cl_{\rho(P)}B$.

\begin{proposition}\label{prop:T}
For all $B,C\subseteq P$, the following statements hold:
\begin{enumerate}[\rm (1)]
\item $T(\varnothing)=\varnothing$ and $\dn B\subseteq T(B)$;
\item if $B\subseteq C$, then $T(B)\subseteq T(C)$;
\item $T(B\cup C)=T(B)\cup T(C)$;
\item $T(B)$ is a lower set.
\end{enumerate}
\end{proposition}

\begin{proof}
Every member of $\Izero(P)$ is nonempty, while each singleton belongs to $\Izero(P)$ and has cut $\dn b$.
This proves (1).
Monotonicity and the lower-set property follow immediately from \eqref{eq:T}.

For (3), monotonicity gives $T(B)\cup T(C)\subseteq T(B\cup C)$.
If $A\in\Izero(P)$ is contained in $B\cup C$, apply Lemma~\ref{lem:partition} in $(P,\upsilon(P))$.
For one of $D=A\cap B$ and $D=A\cap C$, we have $D\in\Izero(P)$ and
\[
 A^\delta=\cl_{\upsilon(P)}A=\cl_{\upsilon(P)}D=D^\delta.
\]
Thus $A^\delta$ is contained in $T(B)$ or $T(C)$.
Taking the union over $A$ proves the reverse inclusion.
\end{proof}

\begin{proposition}\label{prop:fixed-points}
A subset $C\subseteq P$ is $\rho(P)$-closed if and only if $T(C)=C$.
Moreover,
\begin{equation}\label{eq:T-bounds}
 \dn C\subseteq T(C)\subseteq c_P(C)
 \subseteq\cl_{\upsilon(P)}C.
\end{equation}
\end{proposition}

\begin{proof}
By \eqref{eq:rho-closed}, $C$ is $\rho(P)$-closed if and only if $T(C)\subseteq C$.
Since $C\subseteq T(C)$, this is equivalent to equality.
The first inclusion in \eqref{eq:T-bounds} is Proposition~\ref{prop:T}(1).
Since $c_P(C)$ is a fixed point containing $C$, monotonicity gives
$T(C)\subseteq T(c_P(C))=c_P(C)$.
The last inclusion follows from $\upsilon(P)\subseteq\rho(P)$.
\end{proof}

Define the transfinite iterates by
\[
 T^0(C)=C,~ T^{\beta+1}(C)=T(T^\beta(C)),~
 T^\lambda(C)=\bigcup_{\beta<\lambda}T^\beta(C)
\]
for every nonzero limit ordinal $\lambda$.
This sequence is increasing because $T$ is extensive.

\begin{theorem}\label{thm:iteration}
Let $\kappa=\max\{\aleph_0,|P|\}$, with cardinals viewed as initial ordinals.
For every $C\subseteq P$, there is an ordinal $\gamma_C<\kappa^+$ such that
\[
 T^{\gamma_C+1}(C)=T^{\gamma_C}(C)=c_P(C).
\]
In particular, $T^{\kappa^+}=c_P$ as operators on $\Pow(P)$.
\end{theorem}

\begin{proof}
If $T^\beta(C)\subsetneq T^{\beta+1}(C)$ for every $\beta<\kappa^+$, choose
$x_\beta\in T^{\beta+1}(C)\setminus T^\beta(C)$.
The increasing property makes the $x_\beta$ pairwise distinct, giving an injection from $\kappa^+$ into $P$, a contradiction.
Hence the iteration stabilizes at some $\gamma_C<\kappa^+$.
The set $T^{\gamma_C}(C)$ is then a fixed point, so it is $\rho(P)$-closed and contains $C$.

If $F$ is any $\rho(P)$-closed set containing $C$, transfinite induction using monotonicity and $T(F)=F$ gives $T^\beta(C)\subseteq F$ at every stage.
Therefore $T^{\gamma_C}(C)$ is the least closed set containing $C$, namely $c_P(C)$.
Once equality occurs, the sequence remains constant, which gives the uniform equality at $\kappa^+$.
\end{proof}

\begin{corollary}\label{cor:idempotence}
The operator $T$ is idempotent if and only if $T=c_P$.
\end{corollary}

\begin{proof}
If $T^2=T$, then $T(C)$ is a fixed point containing $C$, so
$c_P(C)\subseteq T(C)$.
Together with \eqref{eq:T-bounds}, this gives equality.
The converse follows from idempotence of topological closure.
\end{proof}

The next example shows that idempotence does not hold in general, even for an algebraic complete lattice.

\begin{example}\label{ex:two-steps}
Let $E=\Nat\times\{0,1\}$ and let $P=\Pow(E)$, ordered by inclusion.
For $f\in2^{\Nat}$, put $G_f=\{(n,f(n)):n\in\Nat\}$, and define
\[
\begin{aligned}
 C&=\{G_f:f\text{ is eventually zero}\},\\
 H&=\{B\subseteq E:|B\cap(\{n\}\times\{0,1\})|\leq1
                        \text{ for every }n\in\Nat\}.
\end{aligned}
\]
Then $\rho(P)=\upsilon(P)=\sigma(P)$, but
\begin{equation}\label{eq:two-steps}
 T(C)=\dn C\subsetneq T^2(C)=c_P(C)=H.
\end{equation}
\end{example}

\begin{proof}
The powerset lattice $P$ is algebraic, with finite subsets of $E$ as its compact elements.
For $e\in E$, put $O_e=\{B\subseteq E:e\in B\}$.
These sets generate $\upsilon(P)$, since
\[
 O_e=P\setminus\dn(E\setminus\{e\}),~
 P\setminus\dn A=\bigcup_{e\notin A}O_e.
\]
They also generate $\sigma(P)$: if $B$ belongs to a Scott-open set $U$, the directed family $\Fin(B)$ has supremum $B$, so some finite $F\subseteq B$ belongs to $U$ and
$B\in\{B'\in P:F\subseteq B'\}=\bigcap_{e\in F}O_e\subseteq U$.
Hence $\upsilon(P)=\sigma(P)$, and \eqref{eq:chain} gives $\rho(P)=\upsilon(P)$.

The subspace $C$ is Hausdorff: two distinct graphs differ at some coordinate $n$, and the traces on $C$ of $O_{(n,0)}$ and $O_{(n,1)}$ separate them.
Thus every nonempty $\upsilon(P)$-irreducible subset contained in $C$ is a singleton, giving $T(C)=\dn C$.

The set $H$ is $\upsilon(P)$-closed because
\[
 P\setminus H=\bigcup_{n\in\Nat}\bigl(O_{(n,0)}\cap O_{(n,1)}\bigr).
\]
For every $B\in H$, each finite subset of $B$ extends to an eventually-zero graph.
Consequently $\Fin(B)\subseteq T(C)$.
This is a nonempty directed family with supremum $B$, so $B\in T^2(C)$.
It follows that
\[
 H\subseteq T^2(C)\subseteq c_P(C)\subseteq H.
\]
Finally, the graph of the constant-one sequence belongs to $H$ but is not contained in any eventually-zero graph, and hence does not belong to $T(C)$.
This proves \eqref{eq:two-steps}.
\end{proof}

\begin{proposition}\label{prop:finite-values}
For every finite $F\subseteq P$,
\[
 T(F)=c_P(F)=\dn F.
\]
\end{proposition}

\begin{proof}
The specialization order of $\rho(P)$ is the given order, so $c_P(\{x\})=\dn x$.
Closure preserves finite unions; hence $c_P(F)=\dn F$.
The result now follows from \eqref{eq:T-bounds}.
\end{proof}

View $\Pow(P)$ as a complete lattice under inclusion.
A monotone map on this lattice is Scott continuous exactly when it preserves directed unions.

\begin{theorem}\label{thm:Scott-continuity}
The following conditions are equivalent:
\begin{enumerate}[\rm (1)]
\item $T:\Pow(P)\to\Pow(P)$ is Scott continuous;
\item $c_P:\Pow(P)\to\Pow(P)$ is Scott continuous;
\item $\rho(P)=\alpha(P)$;
\item $T(C)=c_P(C)=\dn C$ for every $C\subseteq P$.
\end{enumerate}
\end{theorem}

\begin{proof}
Suppose (1) holds.
For every $C\subseteq P$, the family $\Fin(C)$ is nonempty and directed under inclusion, with union $C$.
By Proposition~\ref{prop:finite-values},
\[
 T(C)=\bigcup_{F\in\Fin(C)}T(F)
     =\bigcup_{F\in\Fin(C)}\dn F=\dn C.
\]
Thus the fixed points of $T$ are exactly the lower sets.
By Proposition~\ref{prop:fixed-points}, these are the $\rho(P)$-closed sets, so (3) follows.

If (2) holds, the same directed-union argument gives $c_P(C)=\dn C$ for every $C$, which again yields (3).
If (3) holds, then $c_P(C)=\dn C$, and \eqref{eq:T-bounds} yields (4).
Finally, the map $C\mapsto\dn C$ preserves arbitrary unions, so (4) implies both (1) and (2).
\end{proof}

\section{Subset-system derived topologies}\label{sec:subset-systems}

The proof of Proposition~\ref{prop:tau-r-topology} depends on a finite-decomposition property of the testing subsets.
We now isolate this property and apply it to the subset-system construction of Xu, Huang and Yang~\cite{XuHuangYang2026}.
Their Question~7.1 asks which subset systems have an idempotent derived-topology operator.
The results below give a sufficient condition and show that the derived operators for the Rudin and well-filtered determined systems are not idempotent in general.

\subsection{A finite-decomposition condition}

A \emph{subset system} $H$ on $T_0$ spaces assigns to each space $X$ a family of subsets $H(X)$ containing all singletons, such that
\[
 f(A)\in H(Y)~(A\in H(X))
\]
for every continuous map $f:X\to Y$; its action on maps is direct image.
The testing sets in $H(X)$ are taken to be nonempty.
If a subset system also contains the empty set, deleting that member has no effect on the derived topology.
We write $H(X,\tau)$ when the topology needs to be specified.

Following~\cite{XuHuangYang2026}, put
\begin{equation}\label{eq:tau-H}
 \tau_H=\{U\subseteq X:
 U\cap\cl_\tau A\neq\varnothing\Longrightarrow U\cap A\neq\varnothing
 \text{ for every }A\in H(X,\tau)\},
\end{equation}
and let $d_H(\tau)$ be the topology generated by $\tau_H$.
A topology $\tau$ is \emph{$H$-determined} if $d_H(\tau)=\tau$.
The operator $d_H$ is monotone and inflationary, and the least $H$-determined topology containing $\tau$ exists~\cite{XuHuangYang2026}.
With $H$ fixed, we denote this topology by
\[
 \Hd(\tau)=\bigcap\{\theta:\tau\subseteq\theta
 \text{ and }d_H(\theta)=\theta\}.
\]
In particular,
\begin{equation}\label{eq:H-fixed-hull}
 d_H(\tau)=\Hd(\tau)
 ~\Longleftrightarrow~
 d_H(d_H(\tau))=d_H(\tau).
\end{equation}

\begin{definition}\label{def:H-decomposition}
A subset system $H$ has the \emph{finite-decomposition property} if, for every $T_0$ space $(X,\tau)$ and every decomposition
\[
 A=C_1\cup\cdots\cup C_n,
 ~ A\in H(X,\tau),~ n\geq1,
\]
there exists $m\in\{1,\ldots,n\}$ such that
\begin{equation}\label{eq:H-decomposition}
 C_m\in H(X,\tau)
 ~\text{and}~
 \cl_\tau C_m=\cl_\tau A.
\end{equation}
The sets $C_i$ are arbitrary subsets of $X$ and need not themselves belong to $H(X,\tau)$.
\end{definition}

\begin{proposition}\label{prop:H-decomposition}
Suppose that $H$ has the finite-decomposition property, and let $(X,\tau)$ be a $T_0$ space.
Then $\tau_H=d_H(\tau)$ is a topology with the same specialization order as $\tau$.
Its closed sets are exactly the members of
\begin{equation}\label{eq:H-closed-family}
 \mathcal F_H(\tau)=\{F\subseteq X:
 A\in H(X,\tau),\ A\subseteq F
 \Longrightarrow\cl_\tau A\subseteq F\}.
\end{equation}
Moreover, for every $A\in H(X,\tau)$,
\[
 \cl_{\tau_H}A=\cl_\tau A,
 ~ A\in\Irr(X,\tau_H)\subseteq\Irr(X,\tau).
\]
\end{proposition}

\begin{proof}
Taking complements in~\eqref{eq:tau-H} shows that $U\in\tau_H$ if and only if $X\setminus U\in\mathcal F_H(\tau)$.
The family $\mathcal F_H(\tau)$ contains $\varnothing$ and $X$ and is closed under arbitrary intersections.
Let $F_1,F_2\in\mathcal F_H(\tau)$, and suppose that $A\in H(X,\tau)$ is contained in $F_1\cup F_2$.
Apply the finite-decomposition property to
\[
 A=(A\cap F_1)\cup(A\cap F_2).
\]
For some $i\in\{1,2\}$, the set $B=A\cap F_i$ belongs to $H(X,\tau)$ and has the same $\tau$-closure as $A$.
Hence
\begin{equation}\label{eq:H-union-proof}
 \cl_\tau A=\cl_\tau B\subseteq F_i\subseteq F_1\cup F_2.
\end{equation}
This proves closure under binary unions and identifies the closed sets of the topology $\tau_H$.

For $A\in H(X,\tau)$, every $\tau_H$-closed set containing $A$ contains $\cl_\tau A$, by~\eqref{eq:H-closed-family}.
Thus $\cl_\tau A\subseteq\cl_{\tau_H}A$.
The reverse inclusion follows from $\tau\subseteq\tau_H$.
Moreover, if two $\tau_H$-closed sets cover $A$, the argument giving~\eqref{eq:H-union-proof} shows that one of them contains $A$.
Thus $A$ is $\tau_H$-irreducible.

Finally, singleton tests imply that every member of $\tau_H$ is upper for the specialization order of $\tau$.
Consequently,
\[
 \tau\subseteq\tau_H\subseteq\alpha(X,\leq_\tau),
\]
so the specialization orders agree.
The inclusion $\Irr(X,\tau_H)\subseteq\Irr(X,\tau)$ follows from the continuous identity $(X,\tau_H)\to(X,\tau)$.
\end{proof}

Write $\mathcal S$, $\mathcal C$, $\mathcal D$ and $\mathcal R$ for the systems of singletons, nonempty chains, directed subsets and irreducible subsets, respectively.
Chains and directed sets refer to the specialization order.
We also use the Rudin and well-filtered determined subset systems of~\cite{XuShenXiZhao2020}, denoted here by $\RDset$ and $\WDset$.
The symbol $\RDset$ denotes a subset system and is distinct from the abbreviation RD for irreducible determined spaces used in the preceding sections.

For clarity, a nonempty family $\mathcal K$ of nonempty compact saturated subsets of a space is \emph{filtered} if, for any $K_1,K_2\in\mathcal K$, some $K_3\in\mathcal K$ satisfies $K_3\subseteq K_1\cap K_2$.
A $T_0$ space is \emph{well-filtered} if
\[
 \bigcap\mathcal K\subseteq U
 ~\Longrightarrow~
 K\subseteq U\text{ for some }K\in\mathcal K
\]
for every such filtered family and every open set $U$.
A nonempty subset $A$ is a \emph{Rudin set} if $\cl A$ is a minimal closed set meeting every member of some filtered family of nonempty compact saturated sets.
It is \emph{well-filtered determined} if, for every continuous map $f:X\to Y$ into a well-filtered $T_0$ space, there exists a unique $y\in Y$ such that
\[
 \cl_Y f(A)=\cl_Y\{y\}.
\]
We use all subsets with these properties, not only their closed representatives.
The standard inclusions, understood pointwise, are
\begin{equation}\label{eq:H-standard-inclusions}
 \mathcal S\subseteq\mathcal C\subseteq\mathcal D
 \subseteq\RDset\subseteq\WDset\subseteq\mathcal R;
\end{equation}
see~\cite{XuHuangYang2026,XuShenXiZhao2020}.

\begin{proposition}\label{prop:H-standard-decomposition}
Each of the systems $\mathcal S$, $\mathcal C$, $\mathcal D$, $\RDset$, $\WDset$ and $\mathcal R$ has the finite-decomposition property.
Consequently, $\tau_H=d_H(\tau)$ for each of these systems and every $T_0$ topology $\tau$.
\end{proposition}

\begin{proof}
The assertion for $\mathcal S$ is immediate.
Suppose that a directed set has a finite decomposition $A=C_1\cup\cdots\cup C_n$.
Some $C_m$ is cofinal in $A$.
Indeed, otherwise, for each $i$ choose $a_i\in A$ such that no member of $C_i$ lies above $a_i$.
A common upper bound in $A$ of the $a_i$ belongs to some $C_j$, a contradiction.
A cofinal subset of a directed set is directed, and
\[
 C_m\subseteq A\subseteq\dn C_m\subseteq\cl_\tau C_m.
\]
Thus $C_m$ and $A$ have the same closure.
This proves the assertion for $\mathcal D$.
It also proves the assertion for $\mathcal C$, since every subset of a chain is a chain.

The assertion for $\mathcal R$ is Lemma~\ref{lem:partition}.
For either $\RDset$ or $\WDset$, membership depends only on the closure in the fixed ambient topology.
For Rudin sets this follows directly from the definition; for WD sets it follows because continuous images of sets with the same closure have the same closure.
Now apply Lemma~\ref{lem:partition} and~\eqref{eq:H-standard-inclusions}.
The last assertion follows from Proposition~\ref{prop:H-decomposition}.
\end{proof}

\subsection{A sufficient condition for idempotence}

\begin{corollary}\label{cor:H-one-step}
Suppose that $H$ has the finite-decomposition property and $(X,\tau)$ satisfies
\begin{equation}\label{eq:H-test-preservation}
 H(X,\tau_H)=H(X,\tau).
\end{equation}
Then
\[
 \tau_H=d_H(\tau)=\Hd(\tau).
\]
In particular, this equality holds for every $T_0$ topology when
$H\in\{\mathcal S,\mathcal C,\mathcal D,\mathcal R\}$.
\end{corollary}

\begin{proof}
By Proposition~\ref{prop:H-decomposition}, the testing sets in~\eqref{eq:H-test-preservation} have the same closures in $\tau$ and $\tau_H$.
The defining tests for $(\tau_H)_H$ and $\tau_H$ are therefore identical, so $(\tau_H)_H=\tau_H$.
It follows that $\tau_H$ is $H$-determined, and~\eqref{eq:H-fixed-hull} gives the assertion.

For $\mathcal S$, $\mathcal C$ and $\mathcal D$, condition~\eqref{eq:H-test-preservation} follows from preservation of the specialization order.
For $\mathcal R$, it follows from Theorem~\ref{thm:preservation}.
\end{proof}

Condition~\eqref{eq:H-test-preservation} is an additional sufficient hypothesis, not a consequence of the finite-decomposition property.
Functoriality gives only
\begin{equation}\label{eq:H-test-inclusion}
 H(X,\tau_H)\subseteq H(X,\tau),
\end{equation}
by applying $H$ to the continuous identity $(X,\tau_H)\to(X,\tau)$.
Even when every original testing set remains irreducible and retains its closure, it may cease to belong to $H$.
The next example exhibits this for both $\RDset$ and $\WDset$.

\subsection{Two steps for the Rudin and WD systems}

\begin{theorem}\label{thm:H-two-steps}
Let $X=[0,\omega_1)$, where $\omega_1$ is the first uncountable ordinal, and let $\lambda$ be its usual order topology.
Let
\[
 \tau=\{\varnothing\}\cup
 \{X\setminus K:K\text{ is compact in }(X,\lambda)\}
\]
be the co-compact topology, let $\mu$ be the co-countable topology on $X$, and let $\tau_{\mathrm{disc}}=\Pow(X)$ be the discrete topology.
Then $\tau$ is $T_1$, and for each $H\in\{\RDset,\WDset\}$,
\begin{equation}\label{eq:H-two-steps}
 \begin{aligned}
 \tau&\subsetneq\tau_H=d_H(\tau)=\mu,\\
 \mu&\subsetneq\mu_H=d_H(\mu)=\tau_{\mathrm{disc}}=\Hd(\tau).
 \end{aligned}
\end{equation}
In particular, $d_H(\tau)$ need not be $H$-determined, and the iteration of $d_H$, with the original topology at stage $0$, first stabilizes at stage $2$ for both systems.
\end{theorem}

\begin{proof}
Every $\lambda$-compact subset of $X$ is closed and bounded, and hence countable.
Conversely, every countable subset is contained in a compact interval $[0,\beta]$ for some $\beta<\omega_1$.
The compact subsets of the Hausdorff space $(X,\lambda)$ form the proper closed sets of the topology $\tau$.
In particular, $\tau\subseteq\lambda$ and $\tau$ is $T_1$.

We first determine the irreducible subsets of $(X,\tau)$.
For every $\lambda$-compact set $K$, the restrictions of $\tau$ and $\lambda$ to $K$ agree.
Indeed, if $V$ is relatively $\lambda$-open in $K$, then $K\setminus V$ is $\lambda$-compact, and
\[
 V=K\cap\bigl(X\setminus(K\setminus V)\bigr)
\]
is relatively $\tau$-open.
It follows that every countable irreducible subset is a singleton, since it lies in a compact interval on which the induced topology is Hausdorff.
Every uncountable subset is $\tau$-dense, because every proper $\tau$-closed set is countable.
Moreover, any two nonempty $\tau$-open sets intersect, so $X$ is irreducible.
Every dense subset of an irreducible space is irreducible, and hence
\begin{equation}\label{eq:H-omega1-irreducibles}
 \Irr(X,\tau)=\mathcal S(X)\cup
 \{A\subseteq X:A\text{ is uncountable}\}.
\end{equation}

For $\alpha<\omega_1$, set $K_\alpha=[\alpha,\omega_1)$.
Each $K_\alpha$ is $\tau$-compact.
To see this, take a cover by $\tau$-open sets and select a nonempty member $U_0=X\setminus C_0$, where $C_0$ is $\lambda$-compact.
The remainder $K_\alpha\cap C_0$ is $\lambda$-compact because $K_\alpha$ is $\lambda$-closed.
All members of the cover are $\lambda$-open, so finitely many cover this remainder; together with $U_0$ they cover $K_\alpha$.
Since $\tau$ is $T_1$, each $K_\alpha$ is also saturated.
The family $\{K_\alpha:\alpha<\omega_1\}$ is filtered and consists of nonempty sets.
Every proper $\tau$-closed set is bounded and therefore misses some $K_\alpha$.
Consequently, $X$ is the unique closed set meeting every member of this family, and is a Rudin set.

Every uncountable subset has $\tau$-closure $X$, so it too is a Rudin set.
Combining this observation with~\eqref{eq:H-standard-inclusions} and~\eqref{eq:H-omega1-irreducibles} gives
\begin{equation}\label{eq:H-omega1-tests}
 \RDset(X,\tau)=\WDset(X,\tau)
 =\mathcal S(X)\cup\{A\subseteq X:A\text{ is uncountable}\}.
\end{equation}
For a nonempty set $U\subseteq X$, the tests in~\eqref{eq:tau-H} therefore say exactly that $U$ meets every uncountable subset of $X$.
Singleton tests impose no further restriction, because $\tau$ is $T_1$.
Thus
\[
 U\in\tau_H~\Longleftrightarrow~ X\setminus U\text{ is countable},
\]
and $\tau_H=d_H(\tau)=\mu$.
The inclusion $\tau\subseteq\mu$ is strict: the countable set of finite ordinals $\omega$ is not $\lambda$-closed, and hence not $\lambda$-compact, so $X\setminus\omega\in\mu\setminus\tau$.

Every compact subset of $(X,\mu)$ is finite.
Otherwise an infinite compact subset would contain a countably infinite subset $D$, which is closed in $(X,\mu)$ and has the discrete subspace topology.
This would make $D$ an infinite discrete compact space, a contradiction.
It follows that $(X,\mu)$ is well-filtered.
Indeed, in any nonempty filtered family $\mathcal F$ of nonempty finite sets, choose a member $F_0$ of smallest cardinality.
For each $F\in\mathcal F$, there is $G\in\mathcal F$ with $G\subseteq F_0\cap F$.
Minimality of $|F_0|$ forces $G=F_0$, so $F_0\subseteq F$ and $\bigcap\mathcal F=F_0$.
This proves the well-filtered condition.

If $A\in\WDset(X,\mu)$, apply its defining property to the identity map of the well-filtered space $(X,\mu)$.
There is $x\in X$ such that
\[
 \cl_\mu A=\cl_\mu\{x\}=\{x\},
\]
where the last equality uses $T_1$.
Thus $A=\{x\}$, and~\eqref{eq:H-standard-inclusions} yields
\begin{equation}\label{eq:H-cocountable-tests}
 \RDset(X,\mu)=\WDset(X,\mu)=\mathcal S(X).
\end{equation}
Hence every subset of $X$ satisfies the $H$-open tests relative to $\mu$, giving $\mu_H=d_H(\mu)=\tau_{\mathrm{disc}}$.
The inclusion $\mu\subsetneq\tau_{\mathrm{disc}}$ is strict because $X$ is uncountable.

Finally, the discrete topology is $H$-determined.
If $\theta\supseteq\tau$ is any $H$-determined topology, monotonicity gives
\[
 \tau_{\mathrm{disc}}=d_H^2(\tau)\subseteq d_H^2(\theta)=\theta.
\]
Thus $\theta=\tau_{\mathrm{disc}}$ and $\Hd(\tau)=\tau_{\mathrm{disc}}$.
Both preceding inclusions are strict, so the iteration first stabilizes at stage $2$.
\end{proof}

\begin{remark}\label{rem:H-two-steps}
For the topology in Theorem~\ref{thm:H-two-steps}, one also has
\[
 d_{\mathcal R}(\tau)=\mu,
 ~ d_{\mathcal R}(\mu)=\mu.
\]
Indeed, the original Rudin, WD and irreducible testing families coincide by~\eqref{eq:H-omega1-tests}, and Theorem~\ref{thm:one-step} applies to $\mathcal R$.
Thus the passage from $\tau$ to $\mu$ preserves irreducibles and their closures, while the Rudin and WD testing families shrink to the singleton family.
The finite-decomposition property therefore guarantees that $\tau_H$ is a topology but does not guarantee its $H$-determinedness.
Corollary~\ref{cor:H-one-step} and Theorem~\ref{thm:H-two-steps} provide, respectively, a sufficient condition and two standard systems with non-idempotent derived operators for the characterization problem in~\cite[Question~7.1]{XuHuangYang2026}.
\end{remark}

\section{Concluding remarks}\label{sec:conclusion}

The irreducible-derived construction preserves both irreducible subsets and their closures, and therefore produces the least RD refinement in one step.
Applied to the upper topology of a poset, it gives the explicit cut description of the weakest order-compatible RD topology.
The binary-prefix domain supplies the remaining negative answers: saturated subspaces need not preserve either RD or MD.
Together these results resolve the four RD questions and the additional MD question considered here.

The operator $T$ separates two aspects of the construction.
The passage from a topology to its least RD refinement is always idempotent, but a single application of $T$ need not compute the closure of an arbitrary subset.
Example~\ref{ex:two-steps} requires two applications even on an algebraic complete lattice.
The fixed-point description and Theorem~\ref{thm:iteration} give closure in general, while Theorem~\ref{thm:Scott-continuity} identifies precisely when directed unions suffice to determine the operator from its finite values.

For general subset systems, the finite-decomposition condition guarantees that $\tau_H$ is a topology and preserves the closures and irreducibility of the original $H$-sets.
The preservation of their membership in $H$ is a further sufficient condition for one-step determination.
For the Rudin and well-filtered determined systems, the example in Section~\ref{sec:subset-systems} has co-countable first derived topology and discrete second derived topology.
Hence $d_H(\tau)$ need not be $H$-determined, even for these standard irreducible subset systems.
Together with the sufficient condition, this gives partial progress on \cite[Question~7.1]{XuHuangYang2026}.
This concerns iteration of an operator on topologies; the two-step example in Section~\ref{sec:operator} concerns iteration of the cut operator on subsets of a fixed poset.

\section*{Funding}

This work was supported by the National Natural Science Foundation of China (Nos. 12471070, 12071199).

\end{document}